\documentclass[graybox]{svmult}

\usepackage{mathptmx}       
\usepackage{helvet}         
\usepackage{courier}        
\usepackage{type1cm}        
\usepackage{makeidx}         
\usepackage{graphicx}        
\usepackage{multicol}        
\usepackage[bottom]{footmisc}
\usepackage{amsmath, amsfonts, amssymb}
\usepackage{color}
\usepackage{tikz}
\usepackage{wrapfig}
\usepackage{subcaption}
\usepackage{float}

\makeindex             

\begin{document}

\title*{The FVC scheme on unstructured meshes for the two-dimensional Shallow Water Equations}
\titlerunning{FVC scheme on unstructured meshes }
\author{\small Moussa Ziggaf, Mohamed Boubekeur, Imad kissami, Fayssal Benkhaldoun and Imad El Mahi}
\institute{
	Moussa Ziggaf
	\at 
	CSEHS, Mohammed VI Polytechnic University Lot 660, 43150 Bengeurir, Morocco \\
	Universit\'e Sorbonne Paris Nord, LAGA, CNRS, UMR 7539, F-93430, Villetaneuse, France\\
	ENSAO, LM2N, Complexe Universitaire, B.P. 669, 60000 Oujda, Morocco\\
	\email{ Corresponding author (Moussa.ziggaf@um6p.ma)} 
	\and 
	Mohamed Boubekeur
	\at CNRS, UMR 7539, Laboratoire Analyse, G\'eom\'etrie et Applications, LAGA, Universit\'e Sorbonne Paris Nord, F-93430, Villetaneuse, France\\ 
	\email{boubekeur@math.univ-paris13.fr}
	\and 
	Imad kissami \at 
	CSEHS, Mohammed VI Polytechnic University Lot 660, 43150 Bengeurir, Morocco \\ 
	\email{Imad.KISSAMI@um6p.ma}
	\and 
	Fayssal Benkhaldoun\at 
	CSEHS, Mohammed VI Polytechnic University Lot 660, 43150 Bengeurir, Morocco \\
	Universit\'e Sorbonne Paris Nord, LAGA, CNRS, UMR 7539, F-93430, Villetaneuse, France\\
	\email{fayssal@math.univ-paris13.fr}
	\and 
	Imad El Mahi\at 
	CSEHS, Mohammed VI Polytechnic University Lot 660, 43150 Bengeurir, Morocco \\
	ENSAO, LM2N, Complexe Universitaire, B.P. 669, 60000 Oujda, Morocco \\
	\email{Imad.ELMAHI@um6p.ma}
}

\maketitle
\vspace{-2.cm}
\abstract{
The fluid flow transport and hydrodynamic problems often take the form of hyperbolic systems of conservation laws. In this work we will present a new scheme of finite volume methods for solving these evolution equations. It is a family of finite volume Eulerian-Lagrangian methods for the solution of non-linear problems in two space dimensions on unstructured triangular meshes. The proposed approach belongs to the class of predictor-corrector procedures where the numerical fluxes are reconstructed using the method of characteristics, while an Eulerian method is used to discretize the conservation equation in a finite volume framework. The scheme is accurate, conservative and it combines advantages of the modified method of characteristics to accurately solve the non-linear conservation laws with a finite volume method to discretize the equations. The proposed Finite Volume Characteristics (FVC) scheme is also non-oscillatory and avoids the need to solve a Riemann problem. Several test examples will be presented for the shallow water equations. The results will be compared to those obtained with the Roe. 
\keywords{Shallow water equations, Finite volume method, Method of characteristics, FVC scheme, Unstructured meshes. 
\\[2pt]
}
}

\vspace{-.8cm}
\section{Introduction}\label{sec1}
\vspace{-.3cm}
Incompressible Navier-Stokes equations have been widely used in the literature to simulate water
flows including eddy diffusion and Coriolis forces, see for example \cite{codina1999numerical,yan2017free}. However, for free-surface flows these models often become complicated due to the presence of moving boundaries within the  flow domain and also due to the inclusion of hydrostatic pressure. Under certain assumptions these models can be replaced by the well-established shallow water equations. Indeed, the shallow water equations can be derived by depth-averaging the three-dimensional Navier-Stokes equations assuming that the pressure is hydrostatic and the vertical scale is far smaller than the horizontal scale, see \cite{abbott1979elements}. In their depth-averaged form, shallow water equations have been used to model many engineering problems in hydraulics and free-surface  flows including tides in coastal regions, rivers, open channel flows, etc. see for instance \cite{churuksaeva2015mathematical,ozgen2016urban}. Developing highly accurate numerical solvers for shallow water equations presents a challenge due to the non-linear aspect of these equations and their coupling through the source terms. More precisely, the difficulty in these  models lies in the coupling terms involving some derivatives of the physical variables that make the system non-conservative and sometimes non-hyperbolic. A class of Eulerian-Lagrangian methods have also been used in \cite{benkhaldoun2015projection} to solve the two-dimensional shallow water equations. This method avoids the solution of Riemann problem and belongs to the finite volume predictor-corrector type methods. The predictor stage uses the method of characteristics to reconstruct the numerical fluxes whereas the corrector stage recovers the conservation equations in the finite volume framework. Numerical results reported in \cite{benkhaldoun2015projection} for two-dimensional shallow water equations demonstrate that this method is robust and more accurate than the Roe and SRNH schemes, but this previous work  was limited to the \textbf{Cartesian mesh}. In this paper, the method is extended to the \textbf{unstructured mesh}. The results presented here show highly accurate solution by using our proposed finite volume characteristics method and confirm its capability to provide accurate and efficient simulations using unstructured meshes for shallow water flows, including Coriolis forces. This paper is organized as follows. The rotating shallow water equations and their projected speed model are presented in Section \ref{sec2}. In Section \ref{sec3}, the numerical method is formulated for the reconstruction of the FVC scheme. The section \ref{sec4} is devoted to numerical results for several test examples for partial dam-break problem and rotating shallow water equations. Finally, the section \ref{sec5} contains concluding remarks and perspectives.
\vspace{-.5cm}
\section{Mathematical Model}\label{sec2}
\vspace{-.2cm}
\subsection{The rotating shallow water model }
\vspace{-.3cm}
The shallow water equations for the free-surface flow 
in two dimensions with the Coriolis forces are formulated as
\begin{equation}\label{1}
	\left\{ \begin{array}{ll}
		\displaystyle \partial_{t} h + \partial_{x} (hu)+\partial_{y} (hv) = 0\\[0.em]
		\displaystyle \partial_{t} (hu) + \partial_{x} (hu^2+\frac{1}{2}gh^2)+\partial_{y} (huv) =  f_chv\\[0.em]
		\displaystyle \partial_{t} (hv) + \partial_{x} (huv)+\partial_{y} (hv^2+\frac{1}{2}gh^2)= -f_chu\\[0.em]
	\end{array}
	\right.
\end{equation}
where $g$ is the gravitational acceleration, $f_c$ is the Coriolis force, $h$ is the water depth, $u$ and $v$ are the depth-averaged velocities. It is well known that the system (\ref{1}) is strictly hyperbolic with real and distinct eigenvalues. The conservative form of (\ref{1}) is
\begin{equation}\label{2}
		\displaystyle  \partial_{t} W + \nabla \cdot \mathbb{F}(W) = Q(W) 
\end{equation}
 $W= \left( \begin{array}{c} h\\[0.em] hu\\[0.em] hv\\ \end{array}\right)$,\ \ \ 
$\mathbb{F}(W)= \left( \left( \begin{array}{c} hu\\[0.em] hu^2+\frac{1}{2}gh^2\\[0.em] huv\\ \end{array}\right), \left( \begin{array}{c} hv\\[0.em] huv\\ [0.em] hv^2+\frac{1}{2}gh^2\\ \end{array}\right) \right)^T$,\ \ \ 
$Q(W)= \left( \begin{array}{c} 0\\[0.em] f_chv\\[0.em] -f_chu\\ \end{array}\right)$\\
\normalsize  
The system of equations (\ref{2}) has to be solved in a bounded spatial domain $\Omega$, with given boundary and initial conditions.
\vspace{-.5cm}
\subsection{Construction of the projected speed model} 
\vspace{-.3cm}
In this section we adopt the same calculation techniques used in the section 2 of \cite{benkhaldoun2015projection} in order to get the projected speed model. The differential form of the projected speed model is
\begin{equation}\label{4}
	\left\{ \begin{array}{ll}
		\displaystyle \frac{\partial h}{\partial t}   + \frac{\partial  hu{_\eta} }{\partial \eta}  =0,\\[.8em] 
		\displaystyle \frac{\partial hu{_\eta}}{\partial t}  +   \frac{\partial   }{\partial \eta}\left(  hu{_\eta}^2+ \frac{1}{2}gh^2  \right)   = f_chu_{\tau} ,\\[.8em]  
		\displaystyle \frac{\partial hu{_\tau}}{\partial t} +  \frac{\partial   }{\partial \eta}\left(  hu{_\eta}u{_\tau} \right) =-f_chu_{\eta} ,
	\end{array}
	\right.
\end{equation}
The system (\ref{4}) can be rewritten as a transport equation form 
\begin{equation}\label{5}
\displaystyle \frac{\partial \mathbf{U}}{\partial t}(t,X) + u_{\eta}(t,X) \frac{\partial\mathbf{U}}{\partial \eta}(t,X)= \mathbf{F}(\mathbf{U},f_c), \quad  \forall \ \ X =(x,y)\in \Omega \subset \mathbb{R}^2,\ \ t>t_0
\end{equation}
\vspace{-.3cm}\\
with, \ \
$ \displaystyle \mathbf{U}= \left( \begin{array}{c} h\\[0.em] u_{\eta}\\[0.em] u_{\tau}\\ \end{array}\right),\ \ 
\left( \begin{array}{c} u_{\tau}\\[0.em] u_{\eta}\\ \end{array}\right) = \left( \begin{array}{c} vn_x - un_y\\[0.em] un_x +vn_y\\ \end{array}\right),\ \ and \ \ \
\textbf{F}(\mathbf{U},f_c)= \left( \begin{array}{c} -h\partial_{\eta} (u_{\eta}) \\[0.em] -g\partial_{\eta} (h)+f_cu_{\tau}\\[0.em] -f_cu_{\eta}\\ \end{array}\right)$\\
The system of equations (\ref{5}) is used only to reconstruct the numerical fluxes while the finite volume method is applied directly to the conservative system  (\ref{2}), see \cite{sahmim2007sign, benkhaldoun2010new}. 
\vspace{-.55cm}
\section{Finite Volume Characteristics scheme }\label{sec3}
\vspace{-.3cm}
In this section we present the finite volume characteristics method for the
numerical solution of the shallow water equations (\ref{1}). The method consists of two steps and can be interpreted as a predictor-corrector approach. The first step deals with the finite volume discretization of the equations whereas in the second step, the reconstruction of the numerical fluxes is discussed.
\vspace{-.5cm}
\subsection{Finite Volume discretization}
\vspace{-.3cm}
The classical finite volume discretization of the system (\ref{2}) without the bathymetry terms is the volume integral over the total volume of the cell $T_i$, which gives
\begin{equation}\label{6}
	\frac{dW_i}{dt}+\frac{1}{|T_i|}\sum_{j\in N(i)}|\gamma_{ij}|\Phi(W_{ij},\textbf{n}_{ij})=Q_i
\end{equation}
\vspace{-.3cm}\\
where\ \ \  	$\displaystyle W_{i} = \frac{1}{|T_{i}|} \int_{T_{i}} W d V, \quad	\Phi(W_{ij},\textbf{n}_{ij}) \simeq \frac{1}{|\gamma_{ij}|} \int_{\gamma_{ij}} \mathbb{F}(W)\cdot\textbf{n}_{ij} d \sigma,$\\
 $|T_i|$ denotes the area of the cell $T_i$ and $\gamma_{ij}$ is the edge surrounding the cell $T_i$ and
$N(i)$ is the neighbouring triangles of the cell $T_i$. $\Phi(W_{ij},\textbf{n}_{ij})$ is the numerical flux computed at the interface between the cells $i$ and $j$. The intermediate solution $W_{ij}$ is reconstructed using the characteristic method in the predictor stage. The time discretization of (\ref{6}) is performed by a first order explicit Euler scheme. The time domain is divided into $N$ subintervals $[t_n , t_{n+1}]$ with time step $\Delta t  = t_{n+1} - t_n$ for $n = 0, 1,. . . . , N$. $W^n$ is the value of a generic function $W$ at time $t_n$. The fully-discrete formulation of the system (\ref{2}) is given by
\begin{equation}\label{7}
	W^{n+1}_i=W^{n}_i-\frac{\Delta t}{|T_i|}\sum_{j\in N(i)}|\gamma_{ij}|\Phi(W^n_{ij},\textbf{n}_{ij})+\Delta t Q^n_i
\end{equation}
\vspace{-1.3cm}
\subsection{Flux construction}
\vspace{-.3cm}
In the present study, we reconstruct the numerical flux $\Phi(W^n_{ij},\textbf{n}_{ij})$ using the method of characteristics. The fundamental idea of this method is to impose a regular grid at the new time level and to backtrack the flow trajectories to the previous time level, for more details see \cite{roe1986characteristic,seaid2001quasi}. At the previous time level, the quantities that are needed are evaluated by interpolation from their known values on a regular grid.
\vspace{-.5cm}
\subsubsection{Method of characteristics}
\vspace{-.3cm}
The characteristic curves associated with the equation (\ref{5}) are solutions of the initial-value problem
\vspace{-.3cm}
\begin{equation}\label{8}
	\left\{
	\begin{array}{ll}
		\displaystyle \frac{dX^c(t)}{d t} = u_{\eta}( t,X^c(t))\cdot\textbf{n} \quad  t \in [t_n, t_n +\alpha\Delta t],\ \ \alpha > 0 \\[0.4em] 
		
		\displaystyle X^c(t_n + \alpha\Delta t)= X^*
	\end{array}
	\right.
\end{equation}
The solution of (\ref{8}) can be expressed in an integral form as
\begin{equation}\label{9}
	X^c(t_n) = X^* - \int_{t_n}^{t_n + \alpha\Delta t} u_{\eta}(s,X^c(s))\cdot\textbf{n}\ d s
\end{equation}
\vspace{-.4cm}\\
This integral can be calculated using the integral approximation methods. In our simulations we used a first-order Euler method to approximate the integral in (\ref{9}). The numerical fluxes in (\ref{7}) are reconstructed using the solution of the transport equation (\ref{5}) which is given by
\vspace{-.3cm}\\
\begin{equation}\label{10}
	\textbf{U}(t_n+ \alpha\Delta t, X^*) =\textbf{U}(X^c(t_n)) + \int_{t_n}^{t_n + \alpha\Delta t} \textbf{F}(	\textbf{U}(X^c(s),s),f_c)\ d s
\end{equation}
\vspace{-.3cm}\\
where $	\textbf{U}(t_n+ \alpha\Delta t, X^*)$  is the solution at the characteristic feet. It is computed by interpolation of the departure point $X^c(t_n)$ on the mesh. we used the scattered interpolation methods proposed in \cite{amidror2002scattered}. The integral in (\ref{10}) is calculated using the mind-point rule. This approximation is formulated as
\begin{equation}\label{11}
\displaystyle	\textbf{U}^n_{ij} = \hat{\textbf{U}}^n_{ij} + \alpha \Delta t \textbf{F}( \hat{\textbf{U}}^n_{ij},f_c)
\end{equation}
\begin{wrapfigure}[8]{r}{2.9cm}
	\includegraphics[width=2.9cm,height= 2.3cm]{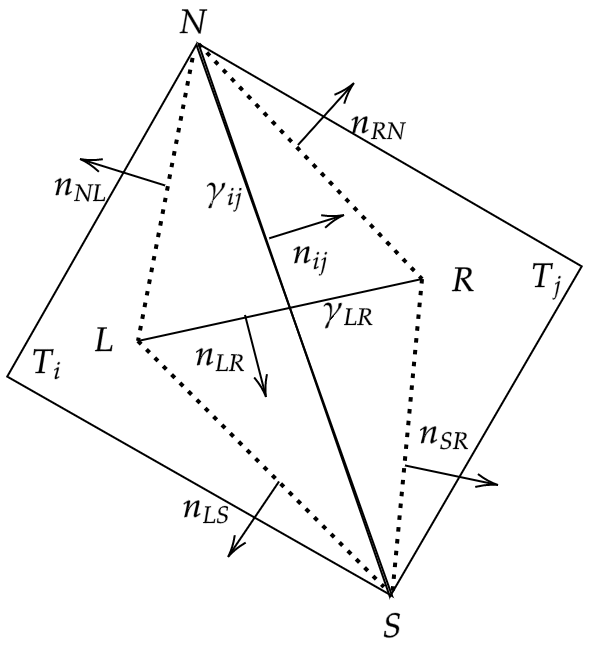}
	\caption{\hspace{-.035cm}Diamond cell.}
 	\label{fig1}
\end{wrapfigure}
where $ \hat{\textbf{U}}^n_{ij}$ is the interpolated solution. To approximate $\textbf{F}(\textbf{U},f_c) $, (i.e. $ \partial_{\eta}(u_{\eta})$,  $ \partial_{\eta}(h),...$ )  we need to approximate these derivatives at the interfaces, for that we use the diamond cell as expressed in Fig \ref{fig1}. For more details see the subsection 3.1.1.2 of \cite{karel2014numerical}. The gradient value at the interface is
\begin{equation}\label{13}
\nabla u_{ij} =\frac{1}{2\mu_{_{SRNL} }} \left\lbrace \frac{}{}(u_{_{S}}-u_{_{N}}) \vec{n_{_{LR}}}|\gamma_{_{LR}}| + (u_{_{R}}-u_{_{L}}) \vec{n_{ij}}|\gamma_{ij}|  \right\rbrace 
\end{equation}\\[-1.em]
where $\mu_{_{SRNL} }$ is the area of the co-volume $SRNL$.
After the discretization of the source term $\textbf{F}(\mathbf{U},f_c)$, the district equations system (\ref{11}) can be written as\\[.3em]
$\displaystyle h^n_{ij} = \hat{h}^n_{ij} - \alpha \Delta t \hat{h}^n_{ij} \nabla  \hat{(u_{\eta} })^n_{ij}$\\[.3em]
$\displaystyle (u_{\eta} )^n_{ij} = \hat{(u_{\eta} })^n_{ij} - \alpha g \Delta t \nabla\hat{h}^n_{ij}+ \alpha \Delta t f_c \hat{(u_{\tau} })^n_{ij}$\\[.3em]
$\displaystyle (u_{\tau} )^n_{ij} = \hat{(u_{\tau} })^n_{ij} - \alpha \Delta t f_c\hat{(u_{\eta} })^n_{ij} $\\[.3em]
Once these projected states are calculated, the states $W_{ij}$
are recovered by using the transformations,
$u^n_{ij} = (u_{\eta} )^n_{ij}n_x -(u_{\tau} )^n_{ij}n_y,\ \ \  \  v^n_{ij} = (u_{\tau} )^n_{ij}n_x + (u_{\eta} )^n_{ij}n_y $\\[.4em]
\textit{ $\blacktriangleright$  The FVC scheme on unstructured meshes for the present model}\\[.4em]
$\left| 
\ \ \begin{aligned}
	\displaystyle W^n_{ij} \  \   \ \  \ \  = \ \ \ & (
	h^n_{ij} \ \ \  h^n_{ij}u^n_{ij} \ \ \  h^n_{ij}v^n_{ij}
	)^T, 
\ \ \ \ \ \  \ \ 
\displaystyle \Phi(W^n_{ij},\textbf{n}_{ij})\ \ =  & \mathbb{F}(W^n_{ij}) \cdot \textbf{n}_{ij} \\[.2em]
\displaystyle W^{n+1}_i \ \  = \ \ \ & W^{n}_i-\frac{\Delta t}{|T_i|}\sum_{j\in N(i)}|\gamma_{ij}|\Phi(W^n_{ij},\textbf{n}_{ij})+\Delta t Q^n_i
\end{aligned}
\right.$

\vspace{-.7cm}
\section{Results} \label{sec4}
\vspace{-.2cm}
In this section we perform numerical tests with our Finite Volume Characteristics scheme on unstructured meshes for the two-dimensional shallow water equations. In all our computations a fixed Courant number $CFL = 0.8$ and $\alpha =1.2 $, are used while the time step
$\Delta t $ is varied according to the stability condition\\
$\displaystyle \Delta t = CFL  \frac{\min_i |\gamma_{ij}|}{\sqrt{2\alpha} \lambda_{ij}^n}, \quad   \lambda_{ij}^n = \max_p \{|u_{pi} ^n+ \sqrt(gh_{pi}^n)|,\  |v_{pj}^n + \sqrt(gh_{pj}^n)|\}$. The used computer is an Intel Core i7-8565U CPU @ 1.80GHz $\times$ 8, with  15  GB RAM.
\vspace{-.7cm}
\subsection{Accuracy test example }
\begin{figure}[H]
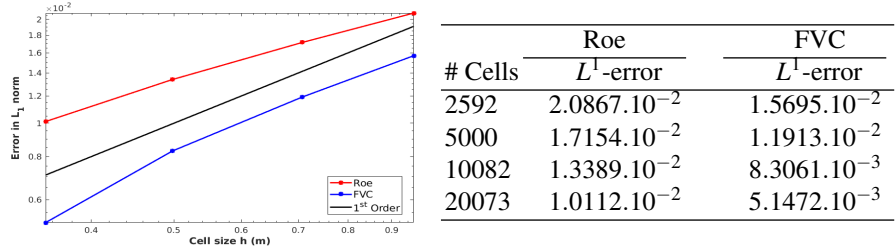

	\caption{ Convergence rates and $L^1$-error. Comparison between FVC and Roe schemes on an unstructured mesh using the same code structure}
\label{fig2}
\end{figure}

\vspace{-.3cm}
The accuracy of the proposed unstructured FVC scheme for a shallow water system is checked, it is compared to the analytical solution. We solve the shallow water equations (\ref{1}) without source terms  in the squared domain $\Omega = [0, 100] \times [0, 100]$ with initial solution for the water depth as the dam-break problem
$h(0,x,y) = 4m,\ \   (x,y) <(0,0) $;\ \  $h(0,x,y) = 2 m, \ \   (x,y) > (0,0) $; \ and \  $ u(0,x,y) = v(0,x,y) = 0 m/s$.
We also compare the results obtained using our FVC scheme on an unstructured mesh to those obtained using the well established Roe scheme in \cite{roe1981approximate}.
The results presented in the table below  are obtained with the relative $L^1$-error norm corresponding to the water depth defined as
$\displaystyle \frac{\sum_{i=1}^{N_{ele}}|T_i| |h^n_i - h(t_n,x_i,y_i)|}{\sum_{i=1}^{N_{ele}}|T_i| |h(t_n,x_i,y_i)|},$
where $h^n_i$ and $h(t_n, x_i , y_i )$ are respectively, the computed and exact water depth at the cell $T_i$, and $N_{ele}$ denotes the total number of cells.
The Relative $L^1$-error is obtained for the accuracy test example at time $t = 5.5 s $ using the Roe and FVC schemes for different unstructured mesh.
We remark that the relative $L^1$-error for the FVC scheme is smaller than for Roe scheme, but the convergence order is still the same and it is close to $1$ (Fig. \ref{fig2}).\\[.5em]
  \begin{minipage}[H]{0.5\linewidth}
	\centering
	\includegraphics[width=6.3cm,height=3.4cm]{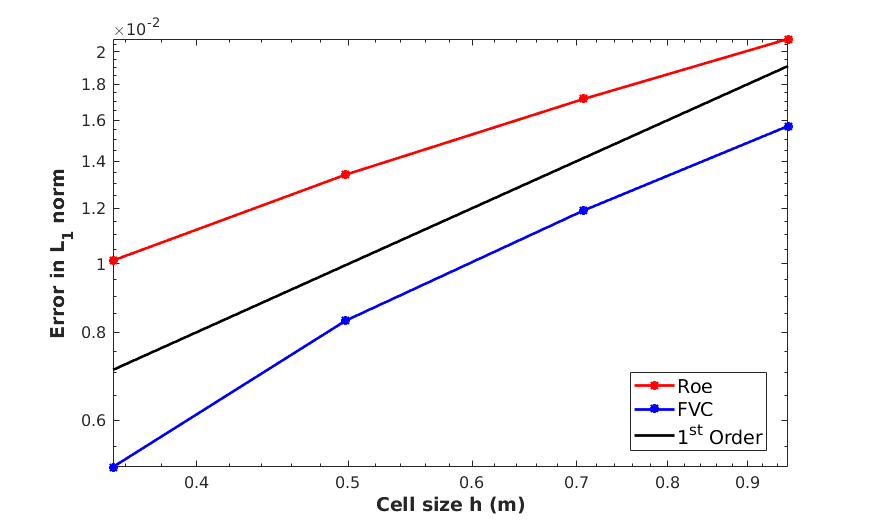}   
\end{minipage}
\begin{minipage}[H]{0.5\linewidth}	
\begin{tabular}{lcccr}
	\hline
	\multicolumn{2}{c}{\hspace{1cm} Roe}&\multicolumn{2}{c}{\hspace{.9cm}FVC}&\\
	\cline{2-2}\cline{4-4}
	$\#$ Cells  & \hspace{.3cm} $L^1$-error& \hspace{.3cm}  \ & \hspace{.3cm} $L^1$-error  & \hspace{.0cm} \\
	\hline
	2592    & \hspace{.2cm}  $2.0867.10^{-2}$ & \hspace{.2cm}\   &  \hspace{.2cm} $1.5695.10^{-2}$  & \\
	5000    &\hspace{.2cm} $1.7154.10^{-2}$ & \hspace{.2cm} \  & \hspace{.2cm} $1.1913.10^{-2}$ & \\
	10082   & \hspace{.2cm} $1.3389.10^{-2}$  & \hspace{.2cm} \  &  \hspace{.2cm} $8.3061.10^{-3}$ & \\
	20073   & \hspace{.2cm}  $1.0112.10^{-2}$  & \hspace{.2cm} \  &  \hspace{.2cm} $5.1472.10^{-3}$ & \\
	\hline
\end{tabular}    
\end{minipage}
\vspace{-.3cm}
\vspace{-1.85cm}
\subsection{Circular dam-break problem }
\vspace{-.2cm}
This benchmark was used in \cite{benkhaldoun2015projection} to represent the FVC scheme on structured Cartesian mesh. We solve the shallow water equations (\ref{1})
on a flat bottom in the spatial domain $\Omega = [-10, 10] \times [-10, 10]$
equipped with the following initial conditions\\
$$h(0,x,y)=1 +\frac{1}{4}\left(1-\tanh \left( \frac{\sqrt{ax^2+by^2}-1}{c}  \right)  \right), \quad  u(0,x,y) = v(0,x,y) = 0 m/s, $$
where $\ a = \frac{5}{2},$\ $b =\frac{2}{5},$\  and\  $c = 0.1$, $g = 1m/s^2 $ and $f_c = 1 Kg.m/s^2 $ as in \cite{benkhaldoun2015projection}. The domain $\Omega$ is discretized with unstructured triangular mesh of 10052 cells.  In this simulation we applied the Neumann conditions on all boundaries (see the subsection 7.5.2 in \cite{mazumder2015numerical}).
\hspace{-12cm}
\begin{figure}[H] 
	\begin{subfigure}[b]{0.33\linewidth}
		\centering
		\includegraphics[width=1\textwidth,height=.12\textheight]{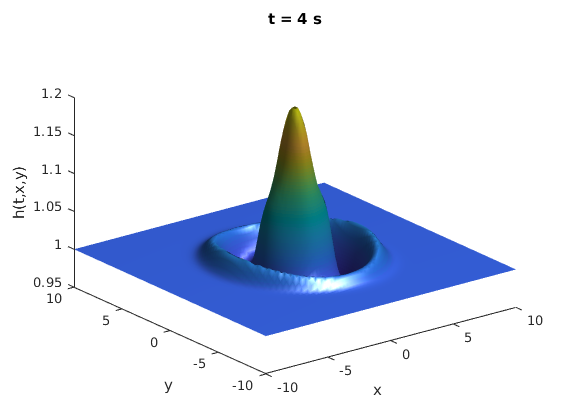} 
	\end{subfigure}
	\begin{subfigure}[b]{0.33\linewidth}
		\centering
		\includegraphics[width=1\textwidth,height=.12\textheight]{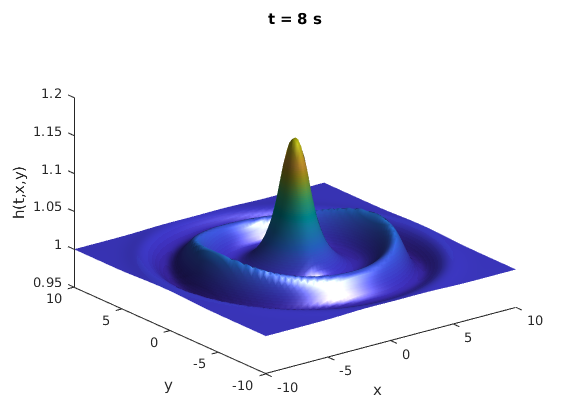} 
	\end{subfigure}
	\begin{subfigure}[b]{0.33\linewidth}
		\centering
		\includegraphics[width=1\textwidth,height=.12\textheight]{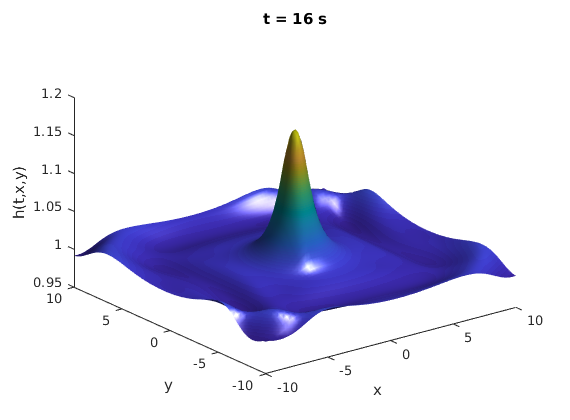} 
	\end{subfigure} 
\end{figure}
\vspace{-1.5cm}
\begin{figure}[H]
	\begin{subfigure}[b]{0.33\linewidth}
		\centering
		\includegraphics[width=1\textwidth,height=.12\textheight]{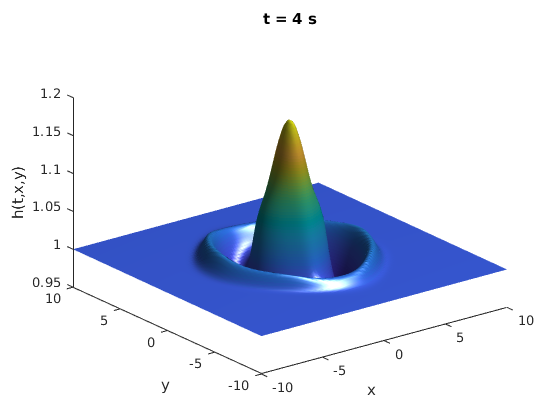}
	\end{subfigure}
	\begin{subfigure}[b]{0.33\linewidth}
		\centering
		\includegraphics[width=1\textwidth,height=.12\textheight]{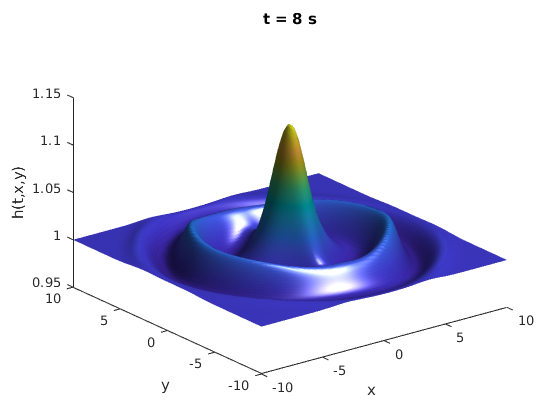}
	\end{subfigure}
	\begin{subfigure}[b]{0.33\linewidth}
		\centering
		\includegraphics[width=1\textwidth,height=.12\textheight]{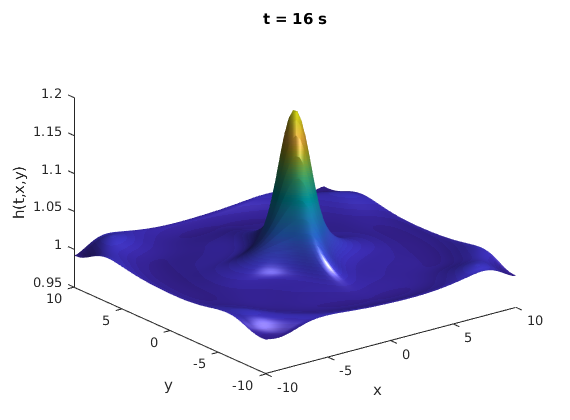} 
	\end{subfigure}
	\caption{Water depth obtained at different times, using FVC on unstructured mesh (first line) and FVC on a Cartesian mesh (second line).}
\label{fig3}
\end{figure}
\vspace{-.6cm}
As it can be clearly seen, the results obtained using FVC scheme on unstructured mesh are very similar to those performed with FVC on structured Cartesian mesh (see subsection 4.2 in \cite{benkhaldoun2015projection}). The rotational movement due to the effect of Coriolis forces
provides an ellipsoid profile, which implies non radial symmetry.  
\vspace{-.6cm}
\subsection{Partial dam-break problem }
\vspace{-.3cm}
This benchmark consists of studying the torrential flow (i.e. Froude number $F_r > 1$) due to a partial and asymmetrical dam-break. This benchmark was proposed in \cite{fennema1989implicit}. Let's study a basin $200 m$ wide, $200 m$ long and flat bottom, without friction.
Water is retained in the left part of the basin.\\[-1em]
\begin{wrapfigure}[13]{r}{3.8cm}
	\vspace{-.8cm}
	\includegraphics[width=4.cm,height=4.cm]{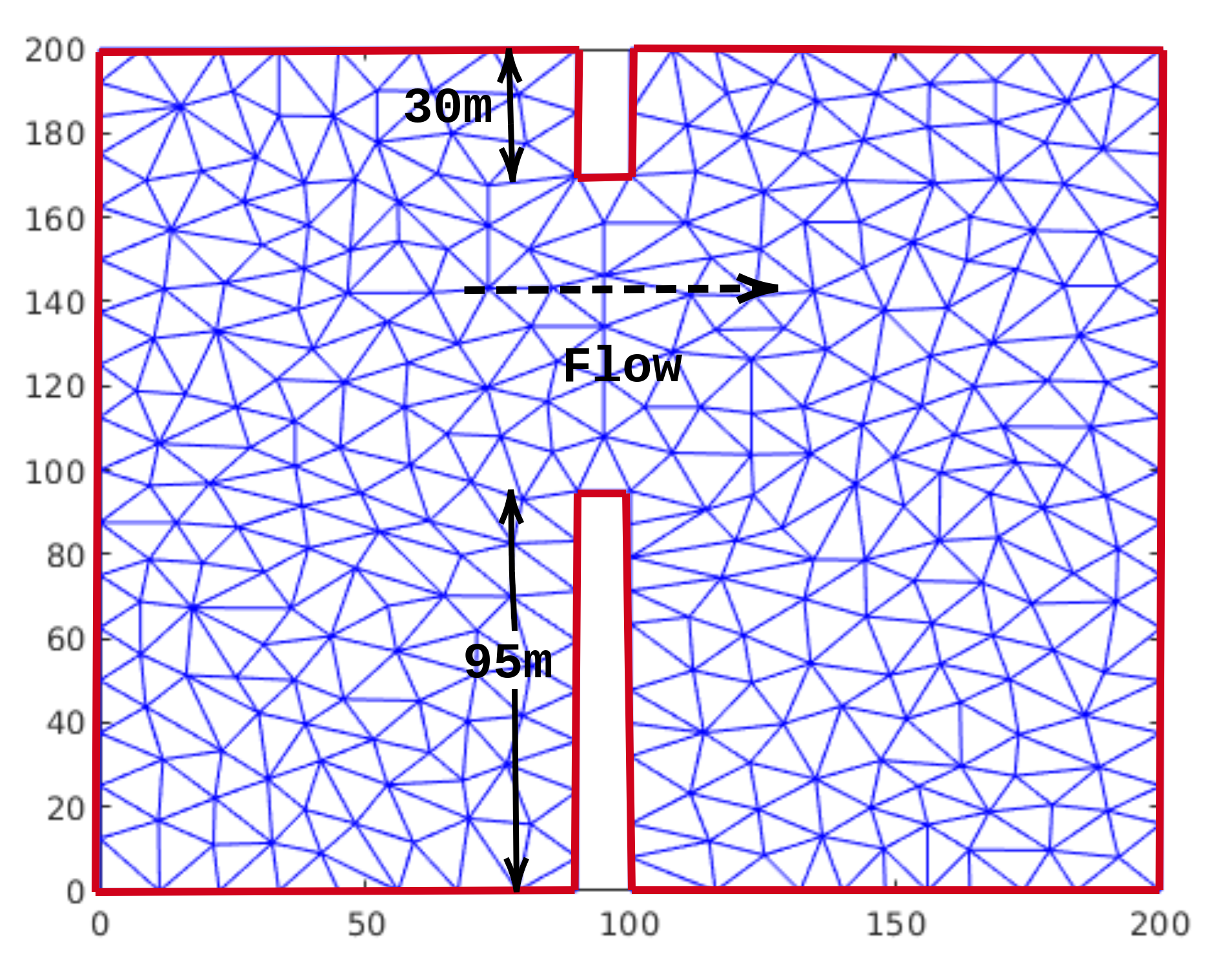}
	\caption{Partial dam-break domain.}
	\label{fig4}
\end{wrapfigure}
The thickness of the dam is  $10 m$ on the flow direction. see Fig \ref{fig4}. Initially  $hr/hl = 0.5$ is fixed with $hl = 4 m$ as water depth
in the reservoir and $hr = 2 m$ as the water level downstream of the dam. The water in the basin is at rest at $t = 0$. When the region occupied by the fluid is bordered by a solid surface, the fluid can not pass through it. Its speed is necessarily zero in the direction perpendicular to the surface. On the other hand, it is not necessarily null in the tangential directions.
In this simulation, a no-slip boundary condition is imposed on all walls see subsection 3.2 in \cite{godlewski2013numerical}. The domain studied was discretized in 20002 non-uniform cells. The duration of simulation is $8.2 s$ counted from the partial dam break.
\vspace{-.5cm}
\begin{figure}[H] 
	\begin{subfigure}[b]{0.33\linewidth}
		\centering
		\includegraphics[width=1\textwidth,height=.15\textheight]{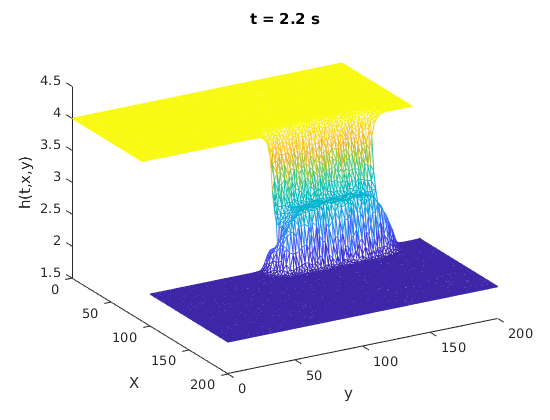} 
	\end{subfigure}
	\begin{subfigure}[b]{0.33\linewidth}
		\centering
		\includegraphics[width=1\textwidth,height=.15\textheight]{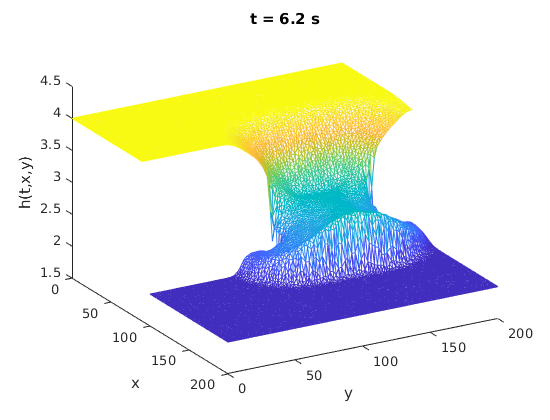} 
	\end{subfigure}
	\begin{subfigure}[b]{0.33\linewidth}
		\centering
		\includegraphics[width=1\textwidth,height=.15\textheight]{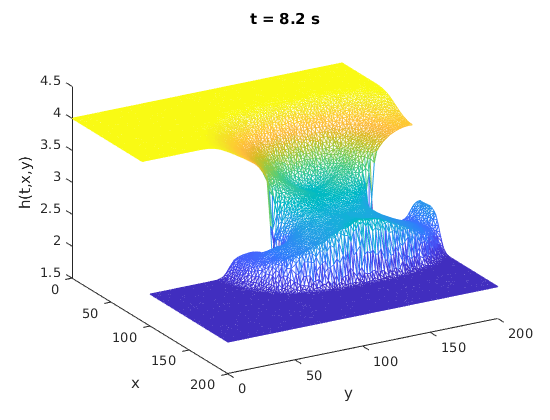} 
	\end{subfigure}
	\caption{Water depth for the partial dam-break problem on flat bottom obtained at different times ($t =2.2s$, $6.2s$ and $8.2s$ ) using FVC scheme on an unstructured mesh. }
\label{fig5}
\vspace{-1.3cm}
\end{figure}
\begin{figure}[H]
	\begin{subfigure}[b]{0.33\linewidth}
		\centering
		\includegraphics[width=1\textwidth,height=.14\textheight]{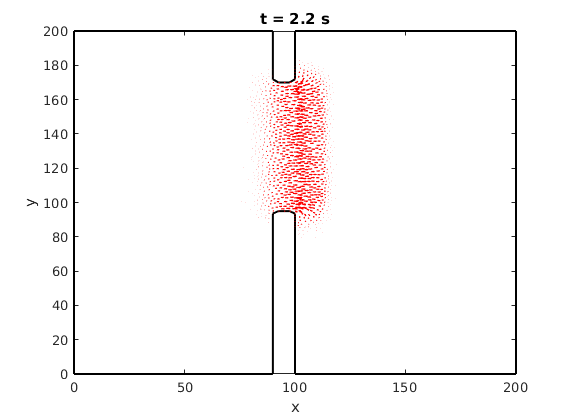} 
	\end{subfigure}
	\begin{subfigure}[b]{0.33\linewidth}
		\centering
		\includegraphics[width=1\textwidth,height=.14\textheight]{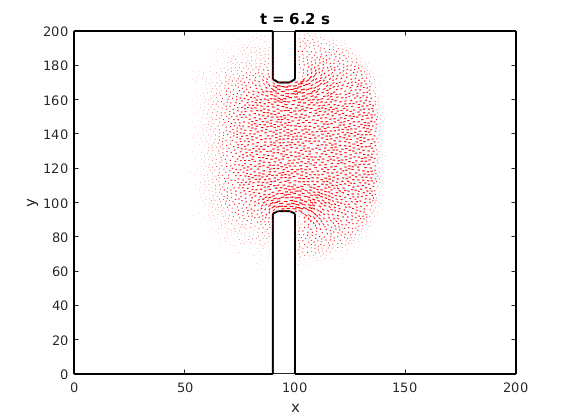} 
	\end{subfigure}
	\begin{subfigure}[b]{0.33\linewidth}
		\centering
		\includegraphics[width=1\textwidth,height=.14\textheight]{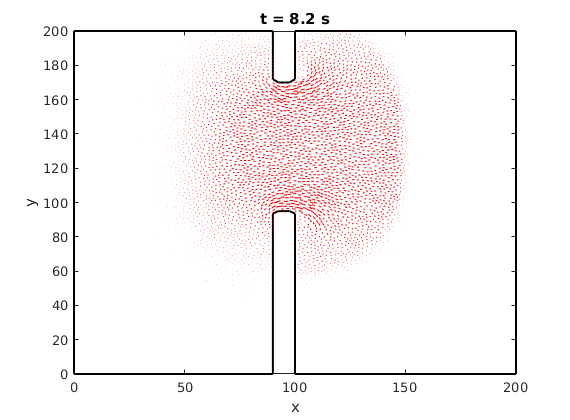} 
	\end{subfigure}
	\begin{subfigure}[b]{0.33\linewidth}
		\centering
		\includegraphics[width=1\textwidth,height=.14\textheight]{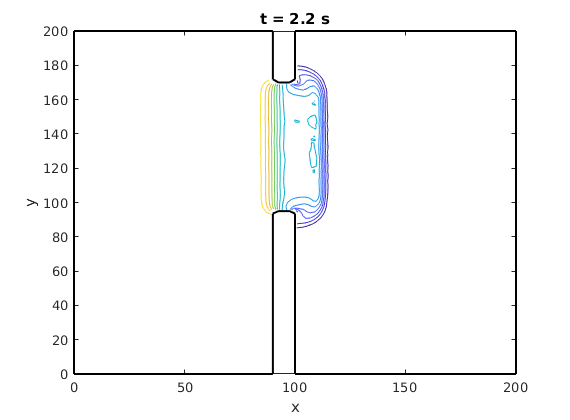} 
	\end{subfigure}
	\begin{subfigure}[b]{0.33\linewidth}
		\centering
		\includegraphics[width=1\textwidth,height=.14\textheight]{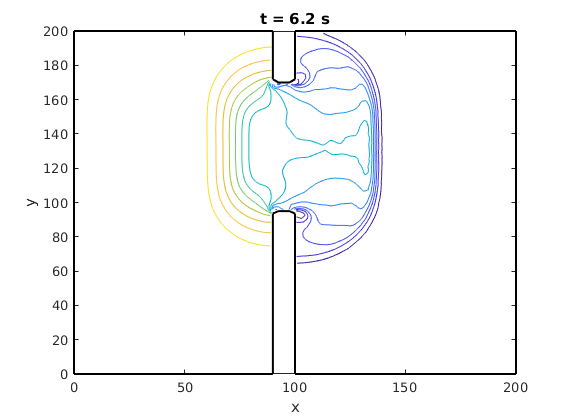} 
	\end{subfigure}
	\begin{subfigure}[b]{0.33\linewidth}
		\centering
		\includegraphics[width=1\textwidth,height=.14\textheight]{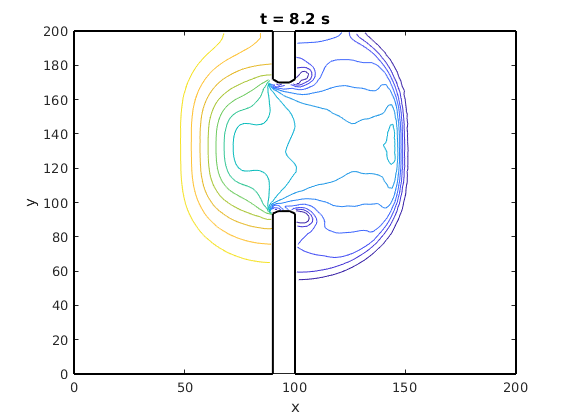} 
	\end{subfigure}
	\caption{Velocity fields and contours for the partial dam-break problem corresponding to the water depth represented in the Fig \ref{fig5}.}
\label{fig6}
\end{figure}

\vspace{-1.cm}
\section{Conclusion} \label{sec5}
\vspace{-.2cm}
A finite volume-characteristics method to solve two-dimensional shallow water equations on unstructured meshes has been presented. This method combines the advantages of the finite volume discretization and the method of characteristics, it solves also steady flows without large numerical errors and compute the numerical flux corresponding to the real state of water flow without relying on Riemann problem solvers. The reasonable accuracy can be obtained easily and no special treatment is needed to maintain a numerical balance, because it is performed automatically in the integrated numerical flux function. Finally, the proposed approach does not require either non-linear solution of algebraic equations or special front tracking techniques. Furthermore, it has strong applicability to various problems in rotating shallow water flows as shown in the numerical results. The outlook of this work is to extend this approach to a multi-layers model of shallow water equations with a bathymetry where we can guarantee a balance between the gradient flux and the source term. In a further step, we will work on coupling this model with the transport convection equation. 

\small
\bibliographystyle{spmpsci}
\bibliography{./Chapitres/references}

\begin{thebibliography}{10}
\providecommand{\url}[1]{{#1}}
\providecommand{\urlprefix}{URL }
\expandafter\ifx\csname urlstyle\endcsname\relax
  \providecommand{\doi}[1]{DOI~\discretionary{}{}{}#1}\else
  \providecommand{\doi}{DOI~\discretionary{}{}{}\begingroup
  \urlstyle{rm}\Url}\fi

\bibitem{abbott1979elements}
Abbott, M.: Elements Of The Theory Of Free Surface Flows; 001.
\newblock Pitman London (1979)

\bibitem{amidror2002scattered}
Amidror, I.: Scattered data interpolation methods for electronic imaging
  systems: a survey.
\newblock Journal of electronic imaging \textbf{11}(ARTICLE), 157--76 (2002)

\bibitem{benkhaldoun2010new}
Benkhaldoun, F., Elmahi, I., Sea{\"\i}d, M.: A new finite volume method for
  flux-gradient and source-term balancing in shallow water equations.
\newblock Computer Methods in Applied Mechanics and Engineering
  \textbf{199}(49-52), 3324--3335 (2010)

\bibitem{benkhaldoun2015projection}
Benkhaldoun, F., Sari, S., Seaid, M.: Projection finite volume method for
  shallow water flows.
\newblock Mathematics and computers in simulation \textbf{118}, 87--101 (2015)

\bibitem{churuksaeva2015mathematical}
Churuksaeva, V., Starchenko, A.: Mathematical modeling of a river stream based
  on a shallow water approach.
\newblock Procedia Computer Science \textbf{66}, 200--209 (2015)

\bibitem{codina1999numerical}
Codina, R.: Numerical solution of the incompressible navier--stokes equations
  with coriolis forces based on the discretization of the total time
  derivative.
\newblock Journal of Computational Physics \textbf{148}(2), 467--496 (1999)

\bibitem{fennema1989implicit}
Fennema, R.J., Hanif~Chaudhry, M.: Implicit methods for two-dimensional
  unsteady free-surface flows.
\newblock Journal of hydraulic research \textbf{27}(3), 321--332 (1989)

\bibitem{godlewski2013numerical}
Godlewski, E., Raviart, P.A.: Numerical approximation of hyperbolic systems of
  conservation laws, vol. 118.
\newblock Springer Science \& Business Media (2013)

\bibitem{karel2014numerical}
Karel, J.: Numerical simulation of streamer propagation on unstructured
  dynamically adapted grids.
\newblock Ph.D. thesis (2014)

\bibitem{mazumder2015numerical}
Mazumder, S.: Numerical methods for partial differential equations: finite
  difference and finite volume methods.
\newblock Academic Press (2015)

\bibitem{ozgen2016urban}
{\"O}zgen, I., Zhao, J., Liang, D., Hinkelmann, R.: Urban flood modeling using
  shallow water equations with depth-dependent anisotropic porosity.
\newblock Journal of Hydrology \textbf{541}, 1165--1184 (2016)

\bibitem{roe1981approximate}
Roe, P.L.: Approximate riemann solvers, parameter vectors, and difference
  schemes.
\newblock Journal of computational physics \textbf{43}(2), 357--372 (1981)

\bibitem{roe1986characteristic}
Roe, P.L.: Characteristic-based schemes for the euler equations.
\newblock Annual review of fluid mechanics \textbf{18}(1), 337--365 (1986)

\bibitem{sahmim2007sign}
Sahmim, S., Benkhaldoun, F., Alcrudo, F.: A sign matrix based scheme for
  non-homogeneous pde’s with an analysis of the convergence stagnation
  phenomenon.
\newblock Journal of Computational Physics \textbf{226}(2), 1753--1783 (2007)

\bibitem{seaid2001quasi}
Sea{\"\i}d, M.: On the quasi-monotone modified method of characteristics for
  transport-diffusion problems with reactive sources.
\newblock Computational Methods in Applied Mathematics Comput. Methods Appl.
  Math. \textbf{2}(2), 186--210 (2001)

\bibitem{yan2017free}
Yan, J., Deng, X., Korobenko, A., Bazilevs, Y.: Free-surface flow modeling and
  simulation of horizontal-axis tidal-stream turbines.
\newblock Computers \& Fluids \textbf{158}, 157--166 (2017)

\end{thebibliography}

\end{document}